\documentclass[twoside,11pt]{amsart}

\usepackage{amsmath} 
\usepackage{amssymb} 

\usepackage{graphicx}
\usepackage{epsfig}
\usepackage{epsf}
\usepackage{float}
\restylefloat{figure}

 \numberwithin{equation}{section} 
 \numberwithin{figure}{section} 
\newtheorem{theorem}{Theorem}[section]

\newtheorem{corollary}[theorem]{Corollary}

\def\be{\begin{equation}}
\def\ee{\end{equation}}
\def\bdm{\begin{displaymath}}
\def\edm{\end{displaymath}}

\def\hf{\frac{1}{2}}

\def\hatf{\widehat{f}}
\def\hatfk{\hatf(k)}

\def\rg{s}
\def\hatgk{\widehat{\rg}(k)}
\def\hatnk{\widehat{n}(k)}

\def\SN{S_N}

\def\RN{R_N}

\def\z{z}  
\def\eps{\varepsilon}
\def\sgn{{\rm sgn}}
\def\gam{k_0}

\def\KNs{K_N^\sigma}

\def\K{{\mathcal K}}
\def\SKNsx{\K_N^\sigma}

\def\lesssim{\stackrel{{}_<}{{}_\sim}}
\def\gtsim{\stackrel{{}_>}{{}_\sim}}
\begin{document}

\title[Spectral data with noise]
{Recovery of edges from spectral data\\ with noise---a new perspective}

\author[Shlomo Engelberg]{Shlomo Engelberg}
\address[Shlomo Engelberg]{\newline
    Department of Electronics\newline 
    Jerusalem College of Technology\newline
    P.O.B. 16031, Jerusalem, Israel}
\email []{shlomoe@jct.ac.il}
\urladdr{http://www.cc.jct.ac.il/\~{}shlomoe/}
\author[Eitan Tadmor]{Eitan Tadmor}
\address[Eitan Tadmor]{\newline
        Department of Mathematics,  Institute for Physical Science and Technology\newline
        and Center of Scientific Computation And Mathematical Modeling (CSCAMM)\newline
        University of Maryland\newline
        College Park, MD 20742 USA}
 \email[]{tadmor@cscamm.umd.edu}
\urladdr{http://www.cscamm.umd.edu/\~{}tadmor}


\date\today

\subjclass{Piecewise smoothness, edge detection, npoisy data, concentration kernels, constrained optimization, separation of scales.}

\thanks{\textbf{Acknowledgment.} E.T. was supported in part by NSF grant 04-07704 and ONR grant N00014-91-J-1076. Part of this research was carried out while the second author was visiting the Weizmann Institute and he is grateful for their hospitality.}


\begin{abstract}
We consider the problem of detecting edges in piecewise smooth
functions from their $N$-degree spectral content, which is assumed to
be corrupted by noise.  There are three scales involved: the
``smoothness" scale of order $1/N$, the noise scale of order $\eta$
and the ${\mathcal O}(1)$ scale of the jump discontinuities.  
We use concentration factors which are adjusted to the noise variance,
$\eta \gg 1/N$, in order to detect the underlying ${\mathcal
O}(1)$-edges, which are separated from the noise scale, $\eta \ll 1$.
\end{abstract}
\maketitle
\tableofcontents

\section{Introduction and statement of main results}\label{Intro}
\setcounter{figure}{0}
\setcounter{equation}{0}

We consider the detection of edges in piecewise smooth data from its Fourier projection,
\[
\SN f(x)=\sum_{|k|\leq N} \hatfk e^{ikx}, \qquad \hatfk:=\frac{1}{2\pi}\int_{-\pi}^\pi f(\xi)e^{-i\xi x}.
\]
Our approach is based on the technique of \emph{concentration
  kernels} advocated in \cite{GeTa99,GeTa00} and the closely related
  techniques described in \cite{Eng}.
The technique presented makes use of the fact that if $f$ is discontinuous, its
Fourier coefficients contain a slowly decaying part associated with
the jumps of $f$. This part decays much more slowly than the rapidly decaying smooth part of $f$.
For example, if $f$ is  smooth except for a single jump discontinuity  at $x=\z$ of size $[f](z):=f(z+)-f(z-)$,
the jump discontinuity is associated with slowly decaying Fourier coefficients,
\[
\hatfk = [f](\z)\frac{e^{-ik\z}}{2\pi ik} + \hatgk, \qquad |\hatgk|
\lesssim \frac{1}{|k|^2}, \quad [f](x):=f(x+)-f(x-), 
\]
where $\hatgk$ are the Fourier coefficients of $\rg$ -- the smooth part of $f$; the smoother $\rg$ is, the faster is the decay of $\hatgk$.  
Concentration kernels succeed in separating the two sets of coefficients.
To this end one computes
\[
\SKNsx * (\SN f)(x)=\pi i\sum_{|k|\leq N} {\rm sgn}(k)\sigma_N\left(\frac{|k|}{N}\right)\hatfk e^{ikx}
\]
Here, $\sigma_N(\xi)$ can be drawn from a large family of properly
normalized concentration factors that are at our disposal. 
The resulting function tends to zero in regions in which $f$ is smooth 
and tends to the amplitude of the jumps at points where the
function has a jump discontinuity,
\[
\SKNsx * (\SN f)(x) = [f](x) + {\mathcal O}\Big(\frac{\log N}{N}\Big). 
\]
Thus, edges are detected by \emph{separation of scales}.

In this paper we utilize concentration kernels to detect edges from spectral information which is corrupted by \emph{white noise}. 
In this context we observe that there are three scales involved -- edges of order ${\mathcal O}(1)$, noise with variance ${\mathcal O}(\eta)$ and the smooth part of $f$ which is resolved within order ${\mathcal O}(1/N)$ or smaller. 
Here, we can separate the noisy part, $n(x)$, from the smooth part,
\[
\rg(x)\rightarrow \rg(x) +n(x), \qquad E(|\hatnk|^2) =\eta;
\]
if $\eta \lesssim 1/N$ then the noisy part could be identified with (or below) the ${\mathcal O}(1/N)$-variation of the smooth part of $f$.
In this case, there are essentially two scales and edges can be detected using the usual framework of concentration kernels advocated in 
\cite{GeTa99,GeTa00}. 
Thus, our main focus in this paper is when the smoothness scale is dominated by the scale of the noise which is still well-separated from  the  ${\mathcal O}(1)$-scale of the jumps,
\[
\frac{1}{N} \lesssim  \eta \ll 1.
\] 
The spectral information is now corrupted by white noise, affecting both low and high frequencies,
\[
\hatfk = [f](\z)\frac{e^{-ik\z}}{2\pi ik} + \hatgk + \hatnk
\]
In order to separate edges from the 
noisy scale, the edge detector, $\KNs f(x)$, must be properly adapted
to the presence of white noise. 
We show how to design edge detectors that
optimally compensate for noise and for the effects of the smooth part of the
signal.

The paper is organized as follows. In section \ref{sec:concentration}
we discuss the general framework of edge detection based on
concentration kernels. We revisit the results of \cite{GeTa00,Ta07},
providing a simpler proof for the concentration property for a large
family of concentration factors. In particular, we trace the precise
dependence of the error on the regularity of the associated
concentration factor.  This will prove useful when we deal with noisy
data in section \ref{sec:perspective}. Here, we introduce our new
perspective, where concentration factors are derived by a constrained
minimization while taking into account the two main ingredients of our
data --- jump discontinuities and the noisy parts of the
data. Numerical results are demonstrated in  
sections \ref{sec:numericsI}. In section \ref{sec:revisited} we extend
our construction of concentration factors to include the \emph{three}
ingredients of the data ---  taking into account the smooth part of
the data in addition to the edges and noisy parts; numerical results
are presented in \ref{sec:numericsII}.  

\section{Detection of edges - concentration kernels}\label{sec:concentration}

Consider an  $f$ which is \emph{piecewise smooth} in the sense that it is sufficiently smooth except for finitely many jump discontinuities, say  at $x=\z_1,\z_2,\ldots,\z_J$, where
\[
[f](\z_j):=f(\z_j+)-f(\z_j-) \neq 0, \quad j=1,2, \ldots J. 
\]
Given the Fourier coefficients, $\big\{\hatfk\big\}_{k=-N}^N$, we are interested in detecting the edges of the underlying piecewise smooth $f$, namely, to detect their location, $\z_1,\ldots, \z_J$ and their amplitudes, $[f](\z_1), \ldots, [f](\z_J)$. 

We utilize edge detection based on \emph{concentration kernels}
\begin{subequations}\label{eqs:KNs}
\begin{equation}\label{eq:ckernel}
\SKNsx(y) := -\sum_{k=1}^N{\sigma_N\left(\frac{k}{N}\right)\sin{ky}}, \qquad \frac{\sigma_N(\xi)}{\xi}\in C^1[0,1].
\end{equation}

We shall need the kernel $-\SKNsx$ to have (approximately) unit mass,
\[
\int_0^\pi \SKNsx(y)dy \approx -1.
\]
To this end, we require that
$\sigma(\xi)\equiv \sigma_N(\xi)$ be a properly normalized
concentration factor; that $\sigma(\xi)$ satisfy
\begin{equation}\label{eq:normalization}
\int_0^1\frac{\sigma(\xi)}{\xi}d\xi=1.
\end{equation}
\end{subequations}
Indeed,  the rectangular quadrature rule yields 
\begin{eqnarray*}
\hspace*{-1cm} \int_0^{\pi}\SKNsx(y)dy&= &\sum_{k=1}^N\sigma\left(\frac{k}{N}\right)\frac{(-1)^k-1}{k}= \nonumber \\
& = & -\sum^N_{ k  \ {\rm odd}\geq 1}\frac{\sigma(\xi_k)}{\xi_k}\frac{2}{N}
=  -\int_0^1\frac{\sigma(\xi)}{\xi}d\xi + \eps_0(N), \quad \xi_k:=\frac{k}{N}. 
\end{eqnarray*}
with an error estimate, e.g., \cite{PZ72,CUN02},
\begin{equation}\label{eq:quad}
\Big|\int_0^{\pi}\SKNsx(y)dy +1\Big| \lesssim \eps_0(N), \qquad \eps_0(N) :=\frac{1}{N}\Big\|\frac{\sigma_N(\xi)}{\xi}\Big\|_{BV} 
\end{equation}
\noindent
We set\footnotemark[1]\footnotetext[1]{We observe that $\KNs f$ is the \emph{operator} associated with but otherwise different from the concentration kernel $\SKNsx(x)$.} 
\[
\KNs f(x):= \SKNsx * f(x) = \pi i\sum_{|k|\leq N} {\rm sgn}(k)\sigma_N\left(\frac{|k|}{N}\right)\hatfk e^{ikx}.
\]
Our purpose is to choose the \emph{concentration factors}, $\sigma_N(|k|/N)$, 
such that $\KNs f(x) \approx [f](x)$. Thus, $\KNs f(x)$ will detect the edges, $[f](\z_j), \ j=1,\ldots, J$, by  concentrating near these ${\mathcal O}(1)$-edges which are to be separated from a  much smaller scale of order $\KNs f(x)\approx 0$ in regions of smoothness. 
In the following theorem we present a rather general framework for edge detectors based on concentration factors. In particular, we track the precise dependence of the scale separation on the behavior of $\sigma$.

\begin{theorem}[Concentration kernels]\label{thm:concentration}
Assume that $f(\cdot)$ is piecewise smooth such that 
the first-variation of $f$ is of locally bounded variation,
\begin{equation}\label{eq:pw}
\omega_f(y)\equiv \omega_f(y;x):=\frac{f(x+y)-f(x-y)-[f](x)}{\sin(y/2)} \in BV[-\pi,\pi].
\end{equation}
Let $\SKNsx(x)$ be an admissible concentration  kernel (\ref{eqs:KNs}), such that
\begin{subequations}\label{eqs:bounds}
\begin{eqnarray}
 \eps_0(N)& :=& \frac{1}{N}\Big\|\frac{\sigma_N(\xi)}{\xi}\Big\|_{BV} \ll 1, \\
 \eps_1(N)&:= &\frac{1}{N}|\sigma_N(1)|  \ll 1,\\
 \eps_2(N)&:=& \frac{1}{N}
\int_{\xi\sim \frac{1}{N}}^1 \frac{|\sigma'_N(\xi)|}{\xi}d\xi  \ll 1 \\
  \eps_3(N)&:= &\Big|\sigma_N\Big(\frac{1}{N}\Big)\Big| \ll 1.
\end{eqnarray}
\end{subequations}
Set 
\begin{equation}
\eps_N:=\eps_0(N)+\eps_1(N)+ \eps_2(N)+ \eps_3(N).
\end{equation}
Then, the conjugate sum $\KNs f(x)$,
\begin{equation}\label{eq:KNs}
\KNs f(x)= \pi i\sum_{|k|\leq N} {\rm sgn}(k)\sigma_N\left(\frac{|k|}{N}\right)\hatfk e^{ikx},
\end{equation}
 satisfies  the concentration property,
\begin{equation}\label{eq:edgedetect}
\KNs  f(x)  \sim  \left\{
\begin{array}{ll} [f](\z_j) + {\mathcal O}\Big( \eps_N\Big),  & x\sim \z_j, \ j=1,\ldots,J\\ \\ 
{\mathcal O}\Big(\eps_N\Big), & dist\Big\{x, \{\z_1, \ldots, \z_J\}\Big\} \gg \eps_N.
\end{array} \right. 
\end{equation}
\end{theorem}

\begin{proof} We simplify the proof in \cite{GeTa00, Ta07}. The key to the proof is
  to observe that $\SKNsx$ is an appropriately normalized
  \emph{derivative} of the delta function; 
in particular, since $\SKNsx(\cdot)$ is odd
\begin{eqnarray*}
\lefteqn{\SKNsx * f(x)  =  -\int_0^{\pi}\KNs (y)\Big(f(x+y)-f(x-y)\Big)dy}   \\
& & =  -\int_0^{\pi}\SKNsx(y) \Big(f(x+y)-f(x-y)-[f](x)\Big)dy
-[f](x)\times \int_0^{\pi}\SKNsx(y)dy. 
\end{eqnarray*}
Since $\sigma$ is assumed normalized, the error estimate (\ref{eq:quad}) tells us\footnotemark[2]\footnotetext[2]{In our case, the errors for the 
trapezoidal and rectangular rules provided respectively in \cite{CUN02} and \cite{PZ72} coincide modulo the ${\mathcal O}\big(\eps_1(N)+\eps_3(N)\big)$ terms.}
\[
 \int_0^{\pi}\SKNsx(y)dy = -1 + \eps_0(N), \qquad \eps_0(N)=\frac{1}{N}\Big\|\frac{\sigma_N(\xi)}{\xi}\Big\|_{BV},
\]
and we end up with the error estimate 
\begin{equation}\label{eq:derr}
\big|\KNs f(x)-[f](x)\big| \lesssim \Big|\int_0^{\pi}\sin(y/2)\SKNsx(y)\omega_f(y;x)dy\Big| + {\mathcal O}\left(\eps_0(N)\right).
\end{equation}
To upperbound the expression on the right, we sum by parts, deriving the identity
\begin{eqnarray*}
\lefteqn{2\sin(y/2)\SKNsx(y) \equiv  \overbrace{\sigma_N(1)\cos{\Big(N +\hf\Big)y}}^{I_1(y)}} \nonumber \\ 
& & + \sum_{k=1}^{N-1}\overbrace{\Big(\sigma_N(\xi_{k}) -\sigma_N(\xi_{k+1})\Big)\cos{\Big(k+\hf\Big)y}}^{I_{2k}(y)} 
\ \  - \overbrace{\sigma_N(\xi_1)\cos\Big(\frac{y}{2}\Big)}^{I_3(y)}.   \nonumber
\end{eqnarray*}
The usual cancelation estimate, ${\displaystyle \big|\int_0^\pi \cos\big((k+\hf)y\big)w(y)dy\big|\lesssim \|w\|_{BV}/|k+\hf|}$, implies 
\begin{eqnarray*}
\Big |\int_0^\pi I_{1}(y)\omega_f(y)dy \Big| & \lesssim & |\sigma_N(1)| \cdot\frac{1}{N}\|\omega_f(\cdot;x)\|_{BV} 
\leq \eps_1(N)\|\omega_f(\cdot;x)\|_{BV}, \\
\Big |\int_0^\pi I_{2k}(y)\omega_f(y)dy \Big| & \lesssim & |\sigma(\xi_{k+1})-\sigma(\xi_k)| \cdot \frac{1}{k+\hf}\|\omega_f(\cdot;x)\|_{BV}, \\
\Big |\int_0^\pi I_3(y)\omega_f(y)dy \Big| & \lesssim & \eps_3(N) \cdot\|\omega_f(\cdot;x)\|_{L^\infty}.
\end{eqnarray*}
We conclude
\begin{eqnarray}\label{eq:ded}
\lefteqn{\Big|\int_0^\pi y\RN(y)\omega_f(y;x)dy\Big|} \nonumber \\
& &   \lesssim \left(\eps_1(N) + \sum_{k=1}^{N-1} \frac{|\sigma(\xi_{k+1})-\sigma(\xi_k)|}{\xi_{k+\hf}}\frac{1}{N} + \eps_3(N)\right)
\|\omega_f(\cdot)\|_{BV}  \\
& & \lesssim \Big( \eps_1(N)+\eps_2(N) + \eps_3(N)\Big)\|\omega_f(\cdot)\|_{BV}. \nonumber
\end{eqnarray}
The result (\ref{eq:edgedetect}) follows from  
(\ref{eq:derr}) and (\ref{eq:ded}).
\end{proof}

As an example, we consider the noise-free case with a concentration
factor for which
\begin{equation}\label{eq:sigmac2}
\frac{\sigma_N(\xi)}{\xi} \in C^1.
\end{equation}
Clearly $\eps_0(N)  \lesssim 1/N$, and since $\sigma$ is bounded, $\eps_1(N) \lesssim 1/N$. Moreover, since $\sigma'_N(\xi)$ is bounded,  
\[
\eps_2(N) \lesssim \frac{1}{N} \int_{\xi\sim \frac{1}{N}}^1
\frac{1}{\xi}d\xi \lesssim \frac{\log(N)}{N}. 
\]
Finally, since $|\sigma_N(\xi)| \lesssim \xi$ then  $\eps_3(N) \lesssim 1/N$. Theorem \ref{thm:concentration} yields the following 
result of \cite{GeTa00, Ta07} for edge detection in piecewise-smooth, noiseless data. 

\begin{corollary}[Concentration kernels]\label{cor:concentration}
Assume that $f(\cdot)$ is piecewise smooth such that (\ref{eq:pw}) holds.
Let $\SKNsx(x)$ be a normalized concentration  kernel (\ref{eqs:KNs}) with $\frac{\sigma(\xi)}{\xi} \in C^1$.
 Then $\KNs f(x)$ satisfies  the concentration property with $\eps_N={\log(N)}/{N}$,
\begin{equation}\label{eq:edgedetectc2}
\KNs  f(x)  \sim  \left\{
\begin{array}{ll} [f](\z_j) + {\mathcal O}\Big( {\displaystyle \frac{\log(N)}{N}} \Big),  & x\sim \z_j, \ j=1,\ldots,J\\ \\ 
{\mathcal O}\Big({\displaystyle \frac{\log(N)}{N}}\Big), & dist\Big\{x, \{\z_1, \ldots, \z_J\}\Big\} \gg {\displaystyle \frac{1}{N}}.
\end{array} \right. 
\end{equation}
\end{corollary}

\section{Noisy data -- a new perspective}\label{sec:perspective}
\subsection{Constrained minimization}
\label{subsect:const_min}
Assume that  $f$ experiences a single jump discontinuity at location
$z$ of height $[f](\z)$.
This dictates a first order decay of the Fourier coefficients,   
\begin{equation}\label{eq:hatnk}
\hatfk = [f](\z)\frac{e^{-ik\z}}{2\pi ik} + \hatgk + \hatnk.
\end{equation}
Here, $\hatgk$ are associated with the regular part of $f$ after
extracting the jump $[f](c)$; their decay is of order  $\sim |k|^{-2}$
or faster, depending on the smoothness of the regular part
$\rg(\cdot)$. The new aspect of the problem enters through the
$\hatnk$'s, which  are the Fourier coefficients of the noisy part
corrupting the smooth part of the data.  We assume $n(\cdot)$ to be
white noise whose mean-square power at each frequency is $E(|\hatnk|^2)=\eta$.
With (\ref{eq:hatnk}), the conjugate sum (\ref{eq:KNs}) becomes
\begin{eqnarray*}
\lefteqn{\KNs f(x)=[f](\z)\times 2\pi i\sum_{k=1}^N \frac{\sigma_N\big(\frac{k}{N}\big)}{k}\cos k(x-\z)} \\
& &  - 2\pi\sum_{k=1}^N \sigma\Big(\frac{k}{N}\Big)\hatgk \sin kx
- 2\pi\sum_{k=1}^N \sigma\Big(\frac{k}{N}\Big)\hatnk\sin kx.
\end{eqnarray*}
 We quantify the ``energy" of each of the three sums on the right. $E_{J}$  
and $E_{R}$ are associated with the discontinuous and regular parts of $f$, 
\begin{subequations}
\begin{eqnarray}
 \sum_{k=1}^N \Big(\frac{\sigma\big(\frac{k}{N}\big)}{k}\Big)^2
\approx \frac{1}{N} \int_0^1 \big(\frac{\sigma(\xi)}{\xi}\Big)^2d\xi=:E_{J}(\sigma),  \label{eq:EJ}\\
\sum_{k=1}^N  \sigma^2\Big(\frac{k}{N}\Big)|\hatgk|^2
\sim \frac{1}{N^3} \int_0^1 \frac{\sigma^2(\xi)}{\xi^4}d\xi =:E_{R}(\sigma),
\label{eq:ER} 
\end{eqnarray}
and $E_{\eta}(\sigma)$ associated with the noisy part of $f$ which was assumed to have variance $\eta$
\begin{equation}\label{eq:EN}
 \sum_{k=1}^N  \sigma^2\Big(\frac{k}{N}\Big)E(|\hatnk|^2) \approx \eta N\int_0^1 \sigma^2(\xi)d\xi =: E_{\eta}(\sigma).
\end{equation}
\end{subequations}

Our perspective for construction of edge detectors for such noisy data
is to treat the problem as a \emph{constrained minimization}. We seek
a function, $\sigma(\xi)$, which minimizes the total energy, thus making the conjugate sum 
$\KNs f$ as localized as possible, subject to  prescribed normalization constraint
(\ref{eq:normalization}),
\begin{equation}\label{eq:minimization}
\min \Big\{a_{J}E_{J}(\sigma)+a_{\eta}E_{\eta}(\sigma)+ a_{R}E_{R}(\sigma) \ \Big| \ \int_0^1\frac{\sigma(\xi)}{\xi}d\xi=1\Big\}.
\end{equation}
This yields
\[
\sigma(\xi)=  C_\sigma \frac{{\displaystyle \frac{1}{\xi}}}{{\displaystyle a_{J}\frac{1}{N\xi^2} +  a_{\eta}\eta N}+ a_{R}\frac{1}{N^3\xi^4}+}= \frac{C_\sigma N^3\xi^3}{a_J N^2\xi^2  + a_\eta\eta N^4\xi^4 + a_R}.
\]
We ignore the relatively negligible contribution of the regular part which becomes even smaller
as $\rg(\cdot)$ becomes smoother. Setting $a_{R}=0$ we end up with concentration factors of the form
\begin{equation}\label{eq:nsigma}
\sigma(\xi) = \frac{C_\sigma}{a_{J}}\cdot\frac{N\xi}{1+ \eta\beta^2 N^2 \xi^2}, \quad 
\beta:=\sqrt{\frac{a_{\eta}}{a_{J}}}.
\end{equation}

The corresponding concentration kernel depends only the \emph{relative} size of the amplitudes $\beta^2=a_{\eta}/a_{J}$: indeed, the normalization
of $\sigma(\xi)/\xi$ (\ref{eq:normalization}) causes the constant
$C_\sigma$ to satisfy
\[
 \int_0^1 \frac{\sigma(\xi)}{\xi}d\xi =
 \frac{C_\sigma}{a_{J}\sqrt{\eta}\beta}\tan^{-1} (\sqrt{\eta}\beta N)
 = 1,
\]
and we end up with the \emph{normalized} concentration factor
\begin{equation}\label{eq:beta}
\sigma(\xi)\equiv \sigma_\eta(\xi)=\frac{1}{\tan^{-1}(\sqrt{\eta}\beta N)}\cdot\frac{\sqrt{\eta}\beta N\xi}{1+\eta\beta^2 N^2\xi^2}, \qquad \beta= \sqrt{\frac{a_{\eta}}{a_{J}}}.
\end{equation}
The corresponding edge detector then takes the form 
\begin{equation}\label{eq:nkernel-a}
K_N^{\sigma_\eta} f(x)=  \frac{\pi \sqrt{\eta}\beta}{\tan^{-1}(\sqrt{\eta}\beta N)} \sum_{|k|\leq N} \frac{ik}{1+\eta\beta^2k^2} \hatfk e^{ikx}.
\end{equation}
The concentration factor $\sigma=\sigma_\eta$ now involves three factors: the ratio $\beta$, the noise variance $\eta$ and the number of modes $N$.
The concentration kernel (\ref{eq:nkernel-a}) tends to de-emphasize both the low frequencies which are ``corrupted" by the jump discontinuity(-ies) and the high frequencies which are corrupted by the noise.
Different procedures yield different policies for the choice of $\beta=\beta(\eta)$; one will be discussed in the next subsection. 
It is worth noting the essential dependence of  $\sigma_\eta(\xi)$ on the variance of the noise $\eta$. There are three scales involved ---  the small ``smoothness" of order $\sim 1/N$, the noise scale of order $\sim \eta$ and the ${\mathcal O}(1)$-scale of jump discontinuities. 
We distinguish between two cases. If $\eta \ll  1/N$ so  that $\sqrt{\eta}\beta N \ll 1$, then the noise can be considered part of the smooth variation of $f$ and 
$\sigma_\eta(\xi) \approx \xi$ recovers the usual concentration factor for noise-free data. Indeed, $\sigma_\eta(\xi)=\xi$ at the limit of $\eta \downarrow 0$. Otherwise, when    the ${\mathcal O}(1/N)$-smoothness scale is dominated by the ${\mathcal O}(\eta)$-noise scale in the sense that $\sqrt{\eta}\beta  \gtsim 1/N$, in which we assume the noise to be still well-below the ${\mathcal O}(1)$-scale of the jumps,
\[
\frac{1}{N} \lesssim \sqrt{\eta}\beta \ll {\mathcal O}(1).
\]
In this case, we can ignore the bounded factor $1/\tan^{-1}(\sqrt{\eta}\beta N)$, and we compute the small scale dictated by theorem 
\ref{thm:concentration}.  Setting $\zeta(\xi):=\sqrt{\eta}\beta N \xi$ we find,
\begin{eqnarray*}
\eps_0(N) & = & \frac{1}{N}\Big\|\frac{\sigma_\eta(\xi)}{\xi}\Big\|_{BV} \lesssim \frac{\sqrt{\eta}\beta N}{N}\int \frac{1}{(1+\zeta^2)^2}d\zeta \lesssim \sqrt{\eta}\beta,  \\ 
\eps_1(N)  & = & \frac{1}{N}|\sigma_\eta(1)|  \lesssim  \frac{1}{\sqrt{\eta}\beta N^2}, \\
\eps_2(N) & =  & \frac{1}{N} \int_{1/N}^1 \frac{1-\zeta^2}{\xi(1+\zeta^2(\xi))^2}
d\xi  \, \lesssim  \sqrt{\eta}\beta\log\big(\sqrt{\eta}\beta\big) \\  
\eps_3(N) & = &  \Big|\sigma_\eta\Big(\frac{1}{N}\Big)\Big|  \lesssim \sqrt{\eta}\beta.
\end{eqnarray*}
and hence (\ref{eqs:bounds}) holds with 
\[
\eps_N\equiv  \eps_\eta:=\sqrt{\eta}\beta \log\big(\sqrt{\eta}\beta\big).
\]
It is remarkable to see how  the small scale of smoothness in the  noiseless case, ${\mathcal O}(\log(N)/N)$,  is now replaced by the small scale of noise 
$\eps_\eta = {\mathcal O}\Big(\sqrt{\eta}\beta \log\big(\sqrt{\eta}\beta\big)\Big)$.
We now appeal to (\ref{eq:edgedetect}): 
since $\eta \ll 1$, theorem \ref{thm:concentration} implies that $K_N^{\sigma_\eta} f$ separates the ${\mathcal O}(1)$-scale of the edges from the noise scale of order $\eps_\eta \ll 1$.

\begin{theorem}[Edge detection in noisy data]\label{thm:noisy}
Assume that $f(\cdot)$ is piecewise smooth in the sense  that (\ref{eq:pw}) holds.
assume that its spectral data contains white noise with variance $\eta \ll 1$.
Let $K_N^{\sigma_\eta} f(x)$ be a normalized concentration  kernel (\ref{eq:nkernel-a}) 
\[
K_N^{\sigma_\eta} f(x)=  \frac{\pi \sqrt{\eta}\beta}{\tan^{-1}(\sqrt{\eta}\beta N)} \sum_{|k|\leq N} \frac{ik}{1+\eta\beta^2k^2} \hatfk e^{ikx}
\]
associated with the concentration factor
\[
\sigma_\eta(\xi)=\frac{1}{\tan^{-1}(\sqrt{\eta}\beta N)}\cdot\frac{\sqrt{\eta}\beta N\xi}{1+\eta\beta^2 N^2\xi^2}.
\]
We distinguish between two cases:\\
(i) if  $\sqrt{\eta}\beta N \ll 1$ we set the small scale $\eps=\eps_N:=\log(N)/N$;\\
(ii) if $1/N \lesssim \sqrt{\eta}\beta  \ll 1$ we set the small scale $\eps=\eps_\eta:=\sqrt{\eta}\beta \log\big(\sqrt{\eta}\beta\big)$.
 Then, $K_N^{\sigma_\eta}f(x)$ satisfies  the following concentration property, 
\begin{equation}\label{eq:edgedetectc3}
\KNs  f(x)  \sim  \left\{
\begin{array}{ll} [f](\z_j) + {\mathcal O}\big( \eps\big),  & x\sim \z_j, \ j=1,\ldots,J\\ \\ 
{\mathcal O}\big(\eps\big), & dist\Big\{x, \{\z_1, \ldots, \z_J\}\Big\} \gg \eps.
\end{array} \right. 
\end{equation}

\end{theorem}

\subsection{Balancing the different types of errors}\label{subsec:balance}

In order to choose the free parameter $\beta$, it is important to know how $\beta$
influences the error at the output of our edge detector.  Let us
consider $E_\eta$.  We have seen that:
\begin{eqnarray*}
E_\eta 
&=& \eta N \int_0^1 \sigma^2(\xi) \, d\xi \\
&=& \frac{\eta N}{(\tan^{-1}(\sqrt{\eta} \beta N))^2}
\int_0^1 
\frac{\eta \beta^2 N^2 \xi^2}{(1 + \eta \beta^2 N^2 \xi^2)^2}
\, d\xi \\
&=&
\frac{\sqrt{\eta}}{\beta (\tan^{-1}(\sqrt{\eta} \beta N))^2}
\int_0^{\sqrt{\eta} \beta N} \frac{\zeta^2}{(1+\zeta^2)^2} \, d\zeta 
\stackrel{\sqrt{\eta} \beta N {\gtsim} 1}{\sim} \frac{\sqrt{\eta}}{\beta}.
\end{eqnarray*}
As $E_\eta$ is approximately the expected value of the square of the
contribution of the noise to the output of the edge detector---it is
approximately the variance of the contribution of the noise, if we
want to consider the {\em size} of the noise, we should consider
something related to the standard deviation of the contribution.  It
is customary to bound the noise by some number of standard
deviations---we will use two standard deviations.  We define the
\emph{effective size} of the noise's contribution to be
\begin{displaymath}
E_{\eta, eff} \equiv 2 \sqrt{E_\eta} \sim  \frac{\eta^{1/4}}{\sqrt{\beta}}.
\end{displaymath}

From the results of Theorem \ref{thm:noisy} we find that far from jump
the contribution from the jump is of order $\sqrt{\eta} \beta$.  The \emph{effective} size of this term is therefore
\begin{displaymath}
E_{J, eff} \sim \eta^{1/2} \beta.
\end{displaymath}

Let us minimize $E_{\eta, eff} + E_{J, eff}$ with respect to $\beta$.
We find that:
\begin{displaymath}
\beta \sim  \eta^{-1/6}.
\end{displaymath}

\section{Numerical results. I}\label{sec:numericsI}
\begin{subequations}\label{figs:noisy}
\begin{figure}[htbp]
\begin{center}
\psfig{file=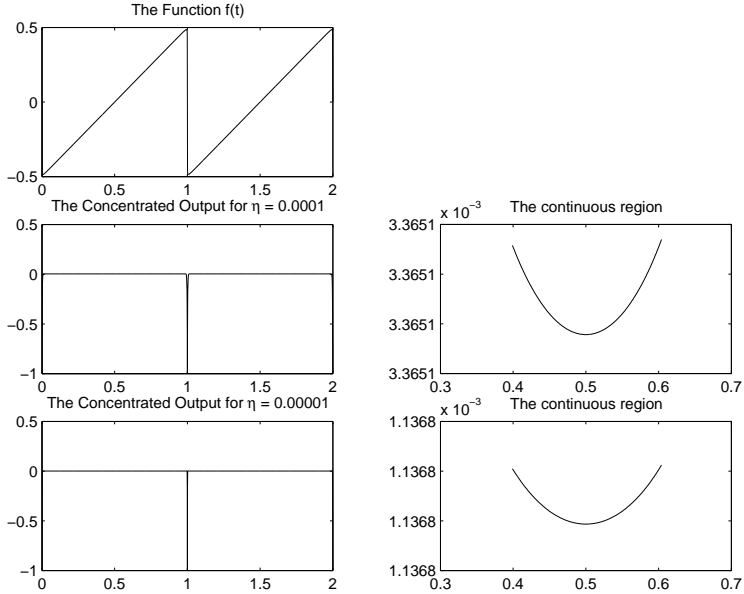,height=8cm}
\end{center}
\caption{The output for various values of $\eta$ when there is
  actually no noise at the input.}
\label{fig:noiseless}
\end{figure}

\begin{figure}[htbp]
\begin{center}
\psfig{file=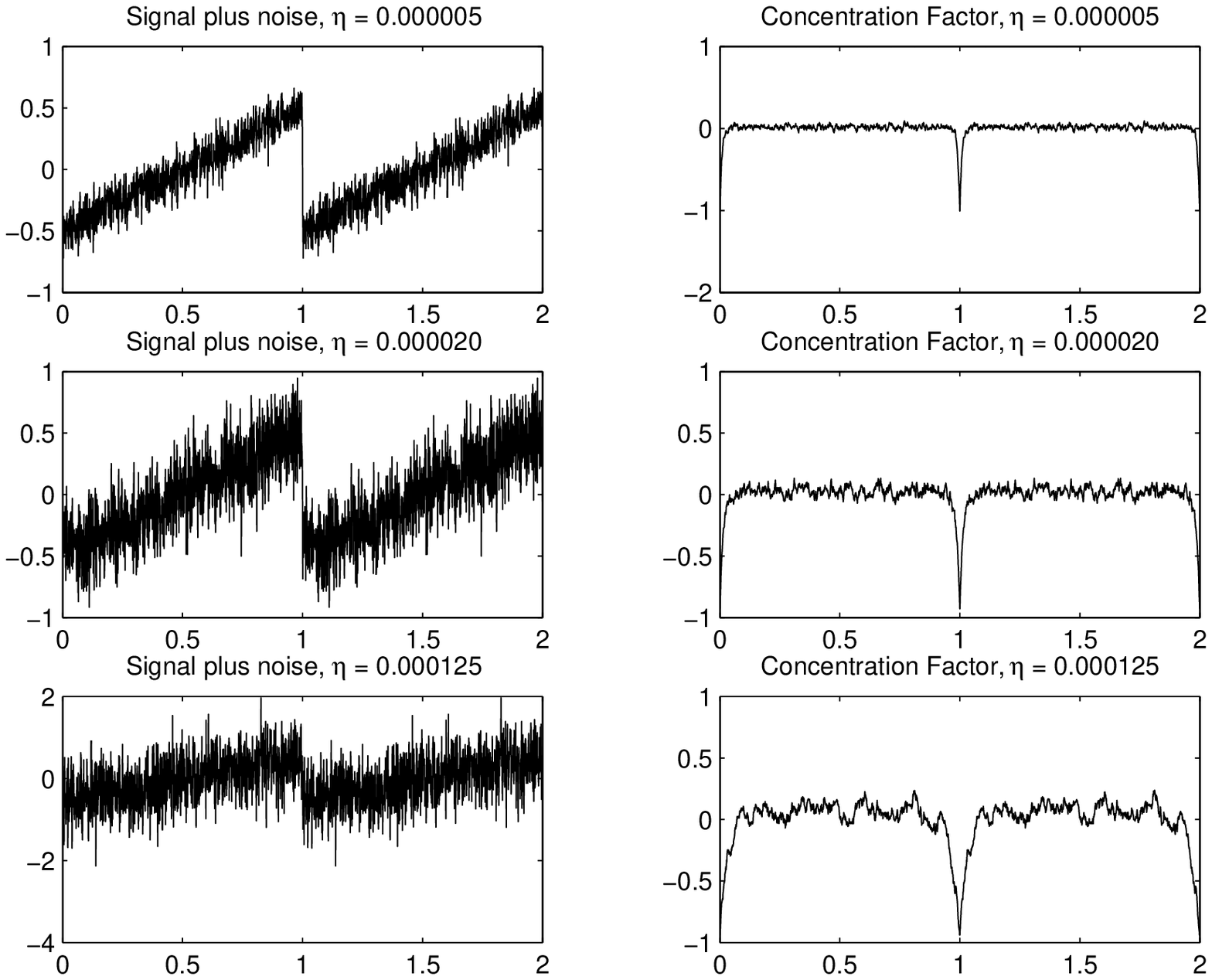,height=10cm}
\end{center}
\caption{Detection of edges in noisy saw-tooth function corrupted with various values of $\eta$, using the concentration kernel (\ref{eq:nkernel-a}) with $\beta=\pi \eta^{-1/6}$.}
\label{fig:noisy-data-a}
\end{figure}
\end{subequations}

To illustrate the results of the previous section, we
present two sets of numerical results.  We begin with the noiseless case, in Figure \ref{fig:noiseless},
where we set $\beta =1$ corresponding to equal weights for the errors
due to the noise and the discontinuous parts of the signal.  We plot
the output of the 
concentration factor when a periodic function with a single jump
continuity is used as the input.  A simple examination of the results
shows that the output is what we predicted.  The output is one at the
(unit) jump and is zero away from it.  As $\eta$ gets smaller the
value away from the jump tends to zero.  Considering the figure, we
find that the ratio of the size in the continuous region is $0.00365 /
0.00137 = 2.66 \approx \sqrt{10}$ which is the ratio of the square root
of the $\eta$'s---as it should be.

Next, Figure \ref{fig:noisy-data-a}
demonstrates the edge detected in noisy data using the concentration kernel (\ref{eq:nkernel-a})
with the advocated $\beta=\pi \eta^{-1/6}$.

\section{Noisy data and smoothness --- concentration kernels revisited}\label{sec:revisited}

As an alternative approach to the $L^2$-minimization offered in Section \ref{sec:perspective}, we now replace the $L^2$-``averaged" effect of the regular part taken in (\ref{eq:ER}), by the BV-like quantity 
\begin{subequations}\label{eq:BV}
\begin{equation}
E_{R}(\sigma)\approx \sum_{k=1}^N \Big|\sigma\Big(\frac{k}{N}\Big)\Big|\cdot|\hatgk|,
\end{equation}
where the regular part is sufficiently smooth that
\begin{equation}
|\hat{s}(k)| \sim \frac{1}{k^2}.
\end{equation}
\end{subequations}
As in (\ref{eq:minimization}), we consider the constrained
minimization
\begin{subequations}\label{eq:2ndmin}
\begin{equation}
\min\Big\{J(\sigma) \  \Big| \ \int_0^1 \frac{\sigma(\xi)}{\xi} \, d\xi = 1 \Big\}, 
\end{equation} 
with  $J(\sigma)=a_J E_J(\sigma)+a_{\eta}E_{\eta}(\sigma)+a_R
E_R(\sigma)$, where $E_J(\sigma)$ and $E_\eta(\sigma)$ given by (\ref{eq:EJ}) and
(\ref{eq:EN}) but with an alternative expression for the ``energy" of
the regular part motivated by (\ref{eq:BV}): $E_{R}(\sigma):= \int_0^1
|\sigma(\xi)|\xi^{-2}d\xi/N$.  We find that:
\begin{equation}
J(\sigma):=  \frac{a_J}{N} \int_0^1 \frac{\sigma^2(\xi)}{\xi^2} \, d\xi
+ a_\eta \eta N \int_0^1 \sigma^2(\xi) \, d\xi
+ \frac{a_R}{N} \int_0^1 \frac{|\sigma(\xi)|}{\xi^2} \, d\xi.
\end{equation}
\end{subequations}
Proceedings formally, the solution for the first variation of (\ref{eq:2ndmin}) leads to
\begin{displaymath}
\frac{2 a_J \sigma(\xi)}{N \xi^2} + 2\eta a_\eta N  \sigma(\xi) +
\frac{a_R \cdot\sgn (\sigma(\xi))}{N \xi^2} = \frac{\lambda}{\xi},\qquad
\sgn(\sigma)=\left\{\begin{array}{rl} 1, & \ \ \sigma >0 \\  0, & \ \
\sigma=0 \\  -1, & \ \ \sigma<0\end{array}\right.  
.
\end{displaymath} 
We will show that the resulting optimal concentration factor is given by
\begin{equation}\label{eq:opt}
\sigma(\xi) = {\displaystyle C_\sigma  \frac{(N \xi - \gam)_+}{1+\eta
    \beta^2 N^2 \xi^2} }, 
\qquad \gam := \frac{a_R}{\lambda}, \ \beta = \sqrt{\frac{a_\eta}{a_J}}, \ 
C_\sigma = \frac{\lambda }{2 a_J}.
\end{equation}
Indeed, to justify the passage to (\ref{eq:opt}), one may consider a
\emph{regularized} 
version of the variational statement (\ref{eq:2ndmin}), $\min J_\eps(\sigma)$, where
\[
J_\eps(\sigma):=  \frac{a_J}{N} \int_0^1 \frac{\sigma^2(\xi)}{\xi^2} \, d\xi
+ a_\eta \eta N \int_0^1 \sigma^2(\xi) \, d\xi
+ \frac{a_R}{N} \int_0^1 \frac{|\sigma(\xi)|_\eps}{\xi^2} \, d\xi,
\]
involves a \emph{mollified} absolute value function:
\begin{displaymath}
|\sigma |_\eps :=
\left\{
\begin{array}{ll}
|\sigma|, & |\sigma| \ge \epsilon, \\
{\displaystyle \frac{\sigma^2}{2 \epsilon} + \frac{\epsilon}{2}}, & |\sigma| \le \epsilon.
\end{array}
\right.
\end{displaymath}
The solution of the corresponding regularized first variation yields the minimizer:
\[
\sigma(\xi) 
=  C_\sigma  \frac{N \xi - \gam \cdot\sgn_\eps(\sigma(\xi))}{\eta \beta^2 N^2 \xi^2 + 1},
\qquad \sgn_\eps(\sigma):=\left\{\begin{array}{rl} 1, & \ \ \sigma >\eps, \\  {\displaystyle {\sigma}/{\eps}}, & \ \ |\sigma|\leq \eps, \\  -1, & \ \ \sigma<-\eps.\end{array}\right.
\]
Thus, we end up with the optimal concentration factor,
$\sigma(\xi)=\sigma_\eps(\xi)$, 
\[
\sigma_\eps(\xi)=\left\{\begin{array}{ll} 
{\displaystyle C_\sigma  \frac{N \xi - \gam}{1+\eta \beta^2 N^2 \xi^2}},  &  N\xi-\gam >\eps, \\ \\
{\displaystyle C_\sigma \frac{\eps N\xi}{\eps(1+\eta \beta^2 N^2 \xi^2) + C_\sigma\gam}}, &
N\xi-\gam\leq \eps, \end{array}\right.
\]
and (\ref{eq:opt}) is recovered by letting $\eps\downarrow 0$.
Clearly, the resulting optimal concentration factor is non-negative. 

It remains to calculate the normalization factor, $C_\sigma$, for which 
\begin{displaymath}
C_\sigma \int_{\gam/N}^1 \frac{1}{\xi} \frac{N \xi - \gam}{\eta \beta^2 N^2   \xi^2 + 1} \, d \xi = 1.
\end{displaymath}
The integral on the left is found to be 
\begin{eqnarray*}
\lefteqn{\int_{{\gam}/{N}}^1 
\left(
\frac{-\gam}{\xi}
+ \frac{\gam \eta \beta^2  N^2  \xi + N}{\eta \beta^2  N^2 \xi^2 + 1}
\right)
\, d\xi=} \\
& & \gam \ln\big(\frac{\gam}{N}\big)
+ \frac{\gam}{2} \ln\left(
\frac{\eta \beta^2  N^2  + 1}{\eta \beta^2 \gam^2+1}
\right)
+
\frac{1}{\sqrt{\eta} \beta} 
\Big(
\tan^{-1}(\sqrt{\eta} \beta N ) - \tan^{-1}(\sqrt{\eta} \beta \gam).
\Big)
\end{eqnarray*}  

We focus our attention on the ``noisy'' case when $1/N \lesssim
\sqrt{\eta} \beta \ll 1$ so that the fourth term on the right is
negligible while the second term on the right is approximated by  
\begin{displaymath}
\frac{\gam}{2} \ln\left(
\frac{\eta \beta^2  N^2  + 1}{\eta \beta^2 \gam^2+1}
\right)
\approx
\gam \ln(\sqrt{\eta}\beta N).
\end{displaymath}  
We end up with an approximated integral
\begin{displaymath}
\int_{{\gam}/{N}}^1 \Big(\ldots \Big)d\xi \approx 
\gam\ln(\sqrt{\eta}\beta N) + \frac{1}{\sqrt{\eta} \beta}\tan^{-1}(\sqrt{\eta} \beta N ).
\end{displaymath}
The balance between these two terms depends on the specific policy for $\beta$ and the detailed balance between $\sqrt{\eta}\beta$ and $N$. Our normalized concentration factor takes the form  

\begin{subequations}\label{eqs:sigeta}
\begin{equation}\label{eq:sigetaa}
\sigma_\eta(\xi)= \frac{1}{\sqrt{\eta}\beta \gam\ln(\sqrt{\eta}\beta N) + \tan^{-1}(\sqrt{\eta} \beta N )}\cdot \frac{\sqrt{\eta}\beta (N\xi-\gam)_+}{1+\eta\beta^2N^2\xi^2}
\end{equation}

\noindent
We can simplify this concentration factor in several ways; we mention
two here.\\ 
(i) When $N$ is large  enough, we have $\tan^{-1}(\sqrt{\eta} \beta N) \approx \pi/2$ yielding
\begin{equation}\label{eq:sigetab}
\sigma_\eta(\xi)= \frac{1}{\gam\ln(\sqrt{\eta}\beta N) + \pi/(2\sqrt{\eta}\beta) }\cdot \frac{(N\xi-\gam)_+}{1+\eta\beta^2N^2\xi^2}
\end{equation}
(ii) Observe that $\sigma_\eta(\xi)$ is rapidly decreasing at $\xi
\approx 1$ with  
\[
\sigma_\eta(\xi) \sim \frac{1}{\eta \beta^2 N \log N}, \qquad \xi \approx 1
\]
so $\sigma_\eta(\xi)$ can be set to zero for $\xi\approx 1$ when $N$
is large enough.  In order to properly normalize the resulting
concentration factor $N$ must be replaced by $N_0$.  This leads us to:
\begin{equation}\label{eq:sigetac}
\sigma_\eta\Big(\frac{k}{N}\Big)= 
\left\{\begin{array}{ll}
0, & k < \gam,\\
{\displaystyle \frac{1}{\gam\ln(\sqrt{\eta}\beta N_0)+\pi/(2\sqrt{\eta}\beta) }\cdot \frac{(k-\gam)_+}{1+\eta\beta^2k^2}}, & 
\gam < k < N_0, \\
0, & N_0 < k <N.
\end{array}\right.
\end{equation}
\end{subequations}

\section{Numerical results. II}\label{sec:numericsII}
\begin{figure}
\begin{center}
\psfig{file=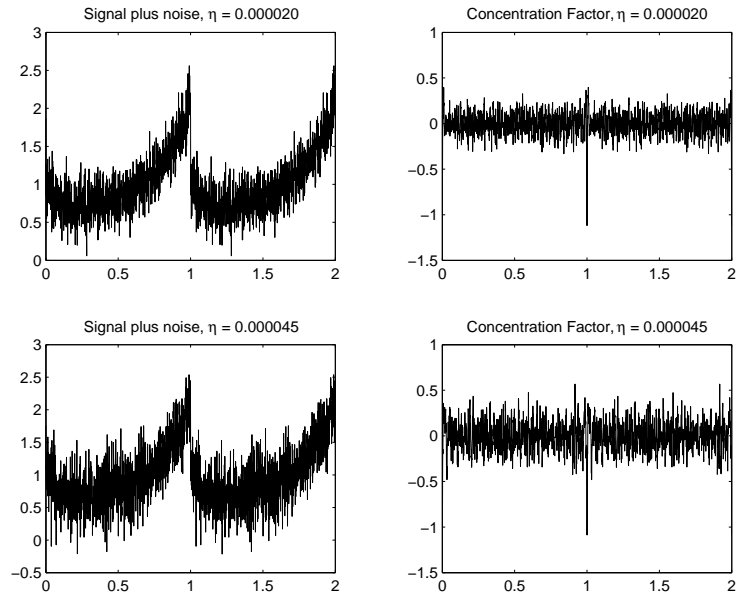,height=8cm}
\end{center}
\caption{The output for various values of $\eta$ when the input
  consists of a piecewise smooth function and noise (of the specified
  power spectral density).} 
\label{fig:smooth_noisy}
\end{figure}

We  consider two  examples depicted in Figure \ref{fig:smooth_noisy}.  
In the first case, we have a noise of variance $\eta = 2\times 10^{-5}$ to be detected out of the first 
$N\gg 1000$ modes. With $\beta=\pi\eta^{-1/6} \sim 15$ and by tuning $\gam=8\pi$ and $N_0=1000$ we find 
\begin{displaymath}
\sigma_\eta\big(\frac{k}{N}\big) =
\left\{
\begin{array}{cl}
0, & k = 1,\ldots,24,1001,1002, \ldots, \\
{\displaystyle 0.3985 \cdot\frac{(k - 8 \pi)_+}{1+\eta(15 k)^2}}, & 
25 \le k \le 1000.
\end{array}
\right.
\end{displaymath}

\noindent
In the second case, of noise variance $\eta = 4.5\times 10^{-5}$ which led us to the choice of  $\beta=\pi\eta^{-1/6} \sim 13$; setting $\gam=6\pi$ and $N_0=1000$ we have
\begin{displaymath}
\sigma_\eta\big(\frac{k}{N}\big)=
\left\{
\begin{array}{cl}
0, & k = 1,\ldots,19, 1001, 1002, \ldots, \\
{\displaystyle 0.5070 \cdot\frac{(k - 6 \pi)_+}{1+\eta (13 k)^2}}, &
20 \le k \le 1000.
\end{array}
\right.
\end{displaymath}
(Note that in calculating the constants we made use of the exact
normalization factor $C_\sigma$.  For our values of $\eta$ and $\beta$
the value $N_0 = 1000$ is not large enough to make the approximate
value given in (\ref{eq:sigetac}) useful.) 

Note that even with a large amount of white
noise and of smooth signal, the location of the jump discontinuity is
still clear.  
When considering jumps ``corrupted'' by low frequency data,
we avoid low frequency signals by not using low frequency
data.  This helps keep the smooth signal from corrupting our results.
On the other hand, because the jump discontinuity has most of its
energy at low frequencies as well, our technique will increase the
noise's effect.  Comparing Figures \ref{fig:noisy-data-a} and
\ref{fig:smooth_noisy}, we find that the latter is not as clean as the
former in the sub-figure where the strength of the noise is the same.

\end{document}